\documentclass[a4,11pt]{article}
\usepackage[latin1]{inputenc}
\usepackage{epsfig}
\usepackage{amsmath}
\usepackage{amssymb}

\newtheorem{theo}{Theorem}[section]

\newtheorem{lema}[theo]{Lemma}

\def\qed{\rule{1.0ex}{1.0ex} \medskip \medskip}

\def \dsp {\displaystyle}

\def\downbar#1{
\setbox10=\hbox{$#1$}
            \dimen10=\ht10 \advance\dimen10 by 2.5pt
            \ifdim \dimen10<15pt %equals approximately 0.5cm
               \advance\dimen10 by -0.5pt
               \dimen11=\dimen10
               \advance\dimen10 by 2.5pt
               \lower \dimen11
            \else \lower \ht10 \fi
            \hbox {\hskip 1.5pt \vrule height \dimen10 depth \dp10}}
\def\upbar#1{
\setbox10=\hbox{$#1$}
            \dimen10=\ht10 \advance\dimen10 by \dp10 \advance\dimen10 by 2.5pt
            \ifdim \dimen10<15pt %equals approximately 0.5cm
               \advance\dimen10 by 2pt \fi
            \raise 2.5pt \hbox {\hskip -1.5pt \vrule height \dimen10}}
\def\cfr#1#2{
\downbar{#2} \hskip -1.2 pt {\; #1 \; \over \thinspace \
#2}\upbar{#1}}
\begin{document}

\title{\Large{ A Favard type theorem for orthogonal polynomials on the unit circle  \\ from a three term recurrence formula}\thanks{This work was support by funds from CAPES, CNPq and FAPESP of Brazil.}}
\author
{
 {K. Castillo$^{a}$, M.S. Costa$^{b}$, A. Sri Ranga$^{a}$\thanks{ranga@ibilce.unesp.br (corresponding author)}\,  and D.O. Veronese$^{c}$} \\[1ex]
  {\small $^{a}$Departamento de Matem\'{a}tica Aplicada, IBILCE, }
  {\small  UNESP - Universidade Estadual Paulista} \\
  {\small 15054-000, São José do Rio Preto, SP, Brazil }\\[1ex]
  {\small $^{b}$Faculdade de Matem\'{a}tica, Universidade Federal de Uberl\^{a}ndia} \\
  {\small  38408-100,  Uberl\^{a}ndia, MG, Brazil }\\[1ex]
    {\small $^{c}$Instituto de Ci\^{e}ncias Tecnol\'{o}gicas e Exatas, Universidade Federal do Tri\^{a}ngulo Mineiro } \\
  {\small  38064-200 Uberaba, MG, Brazil }
}

\date{ }

\maketitle

\thispagestyle{empty}

\begin{abstract}
   The objective of this manuscript is to study directly the Favard type theorem associated with  the three term recurrence formula
\[
    R_{n+1}(z) = \big[(1+ic_{n+1})z+(1-ic_{n+1})\big] R_{n}(z) - 4 d_{n+1} z R_{n-1}(z), \quad n \geq 1,
\]
with $R_{0}(z) =1$ and $R_{1}(z) = (1+ic_{1})z+(1-ic_{1})$, where  $\{c_n\}_{n=1}^{\infty}$ is a real sequence and $\{d_n\}_{n=1}^{\infty}$ is a positive chain sequence.  We establish that there exists an unique nontrivial probability measure $\mu$ on the unit circle for which  $\{R_n(z) - 2(1-m_n)R_{n-1}(z)\}$ gives the sequence of orthogonal polynomials. Here, $\{m_n\}_{n=0}^{\infty}$ is the minimal parameter sequence of the positive chain sequence $\{d_n\}_{n=1}^{\infty}$. The element $d_1$ of the chain sequence, which does not effect the polynomials $R_n$, has an influence in the derived probability measure $\mu$ and hence, in the associated orthogonal polynomials on the unit circle.  To be precise, if $\{M_n\}_{n=0}^{\infty}$ is the maximal parameter sequence of  the chain sequence, then the measure $\mu$ is such that $M_0$ is the size of its mass at $z=1$.  An example is also provided to completely illustrates the results obtained.
\end{abstract}

{\noindent}Keywords: Szeg\H{o} polynomials, Kernel polynomials, Para-orthogonal polynomials, Chain sequences, Continued fractions. \\

{\noindent}2010 Mathematics Subject Classification: 42C05, 33C47.

%\subjclass[2010]{Primary 42C05, 33C47; Secondary 11B83, 11A55, 33C45}

%%%%%%%%%%%%%%%%%%%%%%%%%%%%%%%%%%%%%%%%%%%%%%%%%%%%%%%%%%%%%%%
\setcounter{equation}{0}
\section{Introduction} \label{Sec-Intro}

Orthogonal polynomials on the unit circle (OPUC) have attracted a lot of interest in recent years. For some recent contributions on this topic we refer to  \cite{BreuRyckSimon-2010, Castil-Garza-Marcell-2011, Dam-Mug-Yes-2013, Kheif-Golins-Pehers-Yudit-2011, Pehers-2011, Pehers-Volb-Yudit-2011, Simanek-2012,  Tsujimoto-Zhedanov-2009}. Even though for many years a first hand text for an introduction to these polynomials has been the classical book \cite{Szego-book} of  Szeg\H{o}, detailed accounts regarding the earlier research on these polynomials can be found, for example, in Geronimus \cite{Geronimus-book-1962}, Freud \cite{Freud-book-1971} and Van Assche \cite{VanAssche-Fields-1997}.  However, for recent and more up to date texts on this subject we refer to the two volumes of Simon \cite{Simon-book-p1} and \cite{Simon-book-p2}.  There is also a nice chapter about these polynomials in Ismail \cite{Ismail-book}.

Given a nontrivial  measure $\mu(\zeta) = \mu(e^{i\theta})$ on the unit circle $\mathcal{C} = \{\zeta=e^{i\theta}\!\!: \, 0 \leq \theta \leq 2\pi \}$,  the associated sequence of monic OPUC $\{S_{n}\}$ are those polynomials with the property
\[
\begin{array}l
\dsp \int_{\mathcal{C}} \overline{S_{m}(\zeta)} S_{n}(\zeta) d\mu(\zeta) = \int_{0}^{2\pi} \overline{S_{m}(e^{i\theta})} S_{n}(e^{i\theta}) d\mu(e^{i\theta}) = \delta_{mn}\kappa_n^{-2}.
\end{array}
\]
The orthonormal polynomials on the unit circle are $s_{n}(z) = \kappa_{n} S_{n}(z)$, $n \geq 0$.

The monic  OPUC satisfy the recurrence
\begin{equation} \label{Szego-A-RR}
\begin{array}l
  S_{n}(z) =  z S_{n-1}(z) - \overline{\alpha}_{n-1}\, S_{n-1}^{\ast}(z), \\[1.5ex]
  S_{n}(z) = (1 - |\alpha_{n-1}|^2) z S_{n-1}(z) - \overline{\alpha}_{n-1} S_{n}^{\ast}(z),
\end{array}
n \geq 1,
\end{equation}
where $\overline{\alpha}_{n-1} = - S_{n}(0)$ and $S_{n}^{\ast}(z) = z^{n} \overline{S_{n}(1/\overline{z})}$.   Following Simon \cite{Simon-book-p1} (see also \cite{Simon-2005}) we will be refering to the numbers $\alpha_{n}$ as Verblunsky coefficients.  It is well known that these coefficients  are such that $|\alpha_n| < 1$, $n \geq 0$. It is also known that OPUC are completely characterized in terms of these   coefficients  as given by the following theorem, attributed also to Verblunsky in \cite{Simon-book-p1}. \\

\noindent{\bf Theorem A}
{\em Given an arbitrary sequence  of complex numbers $\{\alpha_{n}\}_{n=0}^{\infty}$, where  $|\alpha_{n}| < 1$, $n \geq 0$, then associated with this sequence there exists a unique nontrivial probability  measure $\mu$ on the unit circle such that  the  monic polynomials $\{S_{n}\}$ generated by $(\ref{Szego-A-RR})$ are the respective monic OPUC. \\
}

In almost all recent studies regarding OPUC, the Verblunsky coefficients and the recurrence relations (\ref{Szego-A-RR}) play a fundamental role. In this manuscript, however, the starting point of the analysis is the  three term recurrence formula
\begin{equation} \label{Eq-TTRR-Rn}
    R_{n+1}(z) = \big[(1+ic_{n+1})z+(1-ic_{n+1})\big] R_{n}(z) - 4 d_{n+1} z R_{n-1}(z), \quad n \geq 1,
\end{equation}
with $R_{0}(z) = 1$ and $R_{1}(z)= (1+ic_{1})z+(1-ic_{1})$, where
\[
  \begin{array}l
     \{c_{n}\}_{n=1}^{\infty} \ \mbox{is a  sequence of real numbers }
  \end{array}
\] \\[-5ex]
\noindent and
\[
  \begin{array}l
\{d_{n}\}_{n=1}^{\infty} \ \mbox{is a  positive chain sequence.}
  \end{array}
\]
For more details on positive chain sequences we refer to Chihara \cite{Chihara-book}.

Although the first element $d_1$ of the chain sequence does not effect the sequence of polynomials $\{R_{n}\}$, its use will become apparent when we introduce the  sequence of polynomials $\{Q_n\}$ in (\ref{Eq-TTRR-Qn}) and, in particular,  the sequence of rational functions $\{A_n/B_n\}$ in section \ref{Sec-Recovering-Szego}.

The main objective of the present manuscript is to show that associated with the sequences $\{c_n\}_{n=1}^{\infty}$ and $\{d_n\}_{n=1}^{\infty}$ there exists an unique nontrivial probability measure $\mu$ on the unit circle (which was also shown in \cite{Costa-Felix-Ranga-2013} by a different method) and to show that the sequence of polynomials $\{R_n(z) - 2(1-m_n)R_{n-1}(z)\}$ are the sequence of OPUC with respect to this measure.  Here, $\{m_n\}_{n=0}^{\infty}$ is the minimal parameter sequence of the positive chain sequence $\{d_n\}_{n=1}^{\infty}$.  As also shown in \cite{Costa-Felix-Ranga-2013}, $M_0$, where $\{M_n\}_{n=0}^{\infty}$ is the maximal parameter sequence of $\{d_n\}_{n=1}^{\infty}$, gives the size of the mass at $z=1$ in the measure $\mu$.

%%%%%%%%%%%%%%%%%%%%%%%%%%%%%%%%%%%%%%%%%%%%%%%%%%%%%%%%%%%%%%%
\setcounter{equation}{0}
\section{Some preliminary results} \label{Sec-PremResults}

The following result gives some information regarding the zeros of $R_{n}(z)$.

\begin{lema} \label{Lemma-SzegoKernel-InterlacingZeros}
The polynomial  $R_{n}(z)$ has  all its $n$ zeros simple and lying on the unit circle $|z| = 1$. Moreover, if one denotes the zeros of $R_{n}$ by $z_{n,j} = e^{i\theta_{n,j}}$, $j=1,2, \ldots,n$,  then
\[
    0 < \theta_{n+1,1} < \theta_{n,1} < \theta_{n+1,2} < \cdots < \theta_{n,n} < \theta_{n+1,n+1} < 2 \pi, \quad n \geq 1.
\]

\end{lema}

This lemma is part of a slightly more extensive result established in \cite{DimRan-2013} with the use of the functions $G_{n}(x)$, defined on the interval  $[-1,1]$, by
\begin{equation} \label{Definition-Gn}
    G_{n}(x) = (4z)^{-n/2} R_{n}(z), \quad n \geq 0,
\end{equation}
where $2x = z^{1/2}+z^{-1/2}$ and $z=e^{i\theta}$. Clearly, the zeros of the function $G_{n}(x)$ in $[-1,1]$ are $x_{n,j} = \cos(\theta_{n,j}/2)$, $j=1,2, \ldots, n$ and Lemma \ref{Lemma-SzegoKernel-InterlacingZeros} means that there holds the interlacing property
\[
   -1 < x_{n+1,n+1} < x_{n,n} < x_{n+1,n} < \cdots < x_{n+1,2} < x_{n,1} < x_{n+1,1} < 1,
\]
for $n \geq 1$. As given in \cite{DimRan-2013}, these functions satisfy the recurrence formula
\begin{equation*} \label{TTRR-for-General-Gn}
    G_{n+1}(x) = \left(x - c_{n+1}\sqrt{1-x^2}\right)G_{n}(x) - d_{n+1}\,G_{n-1}(x), \quad n \geq 1,
\end{equation*}
with $G_{0}(x) = 1$ and $G_{1}(x) = x - c_{1}\sqrt{1-x^2}$.  Moreover, the associated Christoffel-Darboux functions or Wronskians
\begin{equation} \label{Wronskian-for-General-Gn}
   W_{n}(x) = G_{n}^{\prime}(x)G_{n-1}(x)-G_{n-1}^{\prime}(x)G_{n}(x), \quad n \geq 1,
\end{equation}
which not necessarily remain positive through out $[-1,1]$, but satisfy at the zeros of $G_n(x)$
\begin{equation} \label{Wronskian-for-General-Gn-Positiveness}
     W_{n}(x_{n,j}) >  0, \quad W_{n+1}(x_{n,j}) = d_{n+1} W_{n}(x_{n,j}) >  0, \quad j=1,2, \ldots, n,
\end{equation}
for $n \geq 1$. Note that
\[
 W_{n}(x_{n,j}) = G_{n}^{\prime}(x_{n,j})G_{n-1}(x_{n,j}) \quad \mbox{and} \quad W_{n+1}(x_{n,j}) = -G_{n}^{\prime}(x_{n,j})G_{n+1}(x_{n,j}).
\]
Now  we consider the Wronskians
\begin{equation} \label{Wronskian-for-General-Rn}
      V_{n}(z) = R_{n}^{\prime}(z) R_{n-1}(z) - R_{n-1}^{\prime}(z) R_{n}(z), \quad n \geq 1,
\end{equation}
associated with the polynomials $R_{n}(z)$.  From (\ref{Definition-Gn})
\[
     G_{n}^{\prime}(x) = (4z)^{-(n-1)/2} \left[2z R_{n}^{\prime}(z) - n R_{n}(z)\right] \frac{1}{z-1}.
\]
Hence, $ W_{n}(x) = \frac{(4z)^{-(n-1)}}{z-1}\left[2zV_{n}(z) -R_{n-1}(z)R_{n}(z)\right]$, \ $n \geq 1$ \ and
\begin{equation} \label{Wronskian-identity}
 \frac{z_{n,j}^{-(n-2)}}{z_{n,j}-1}V_{n}(z_{n,j}) = 2^{2n-3}W_{n}(x_{n,j})  , \quad j = 1, 2, \ldots, n, \quad n \geq 1.
\end{equation}

From the recurrence formula for $\{R_{n}(z)\}$,
\[
   \frac{R_{n}(1)}{2R_{n-1}(1)} \big[1 - \frac{R_{n+1}(1)}{2R_{n}(1)}\big] = d_{n+1}, \quad n \geq 1.
\]
Hence, $\{\hat{m}_{n}\}_{n=0}^{\infty}$, with
\[  \hat{m}_{n} = 1 - \frac{R_{{n+1}}(1)}{2R_{n}(1)}, \quad n \geq 0, \]
is the minimal parameter sequence of the chain sequence $\{d_{1,n}\}_{n=1}^{\infty}$, where $d_{1,n} = d_{n+1}$, $n \geq 1$.

If we denote by $\{m_n\}_{n=0}^{\infty}$ and $\{M_n\}_{n=0}^{\infty}$ the minimal and maximal parameter sequences of $\{d_n\}_{n=1}^{\infty}$, respectively, then with $m_{1,n} = m_{n+1}$, $M_{1,n} = M_{n+1}$, $n \geq 0$, the sequences $\{m_{1,n}\}_{n=0}^{\infty}$ and $\{M_{1,n}\}_{n=0}^{\infty}$ are parameter sequences of $\{d_{1,n}\}_{n=1}^{\infty}$. Clearly, $\{m_{1,n}\}_{n=0}^{\infty}$ is such that $\hat{m}_n <  m_{1,n}$, $n \geq 0$.  However,
$\{M_{1,n}\}_{n=0}^{\infty}$  is exactly the maximal parameter sequence of $\{d_{1,n}\}_{n=1}^{\infty}$. These and other interesting results on positive chain sequences  can be found in \cite{Chihara-book}.

Note that, the chain sequence $\{d_n\}_{n=1}^{\infty}$ can be such that  $M_0 = m_0 = 0$. But, it is important to note that always $0 < m_{1,0} \leq M_{1,0} < 1$. The equality $m_{1,0} = M_{1,0}$ holds when the chain sequence $\{d_n\}$ has an unique parameter sequence.

%%%%%%%%%%%%%%%%%%%%%%%%%%%%%%%%%%%%%%%%%%%%%%%%%%%%%%%%%%%%%%%
\setcounter{equation}{0}
\section{Some correspondence and asymptotic properties} \label{Sec-CorrAsymProps}

Let $\{Q_n\}$ be the sequence of polynomials given by the continued fraction expression
\[
  \begin{array} {l}
    \dsp \frac{Q_{n}(z)}{R_{n}(z)}  = \cfr{2d_1}{(1+ic_1)z+(1-ic_1)}-\cfr{4d_2z}{(1+ic_2)z+(1-ic_2)}   -  \cdots - \cfr{4d_{n}z}{(1+ic_n)z+(1-ic_n)} \, .
  \end{array}
\]
From the theory of continued fractions (see, for example, \cite{CuytETC-book, JT-Book, LW-book}) the polynomials  $Q_n$ are such that
\begin{equation} \label{Eq-TTRR-Qn}
    Q_{n+1}(z) = \big[(1+ic_{n+1})z+(1-ic_{n+1})\big] Q_{n}(z) - 4 d_{n+1} z Q_{n-1}(z), \quad n \geq 1,
\end{equation}
with $Q_{0}(z) = 0$ and $Q_{1}(z) = 2d_1$.

First we look at an asymptotic result associated with the sequence $\{Q_{n}(1)/R_{n}(1)\}$. From the recurrence formulas for $\{R_{n}(z)\}$ and $\{Q_{n}(z)\}$ together with the theory of continued fractions
\[
  \begin{array} {ll}
    \dsp \frac{Q_{n}(1)}{R_{n}(1)} & \dsp = \cfr{d_1}{1}-\cfr{d_2}{1}  - \cfr{d_3}{1} -  \cdots - \cfr{d_{n}}{1} \, , \\[2ex]
    & \dsp = (1-M_0)\, \cfr{M_{1,0}}{1}-\cfr{(1-M_{1,0})M_{1,1}}{1}  -   \cdots - \cfr{(1-M_{1,n-2})M_{1,n-1}}{1}\, ,
  \end{array}
\]
for $n \geq 1$. Hence, one can write (see the proof of Lemma 3.2 in \cite{Chihara-book})
\[
  \begin{array} {ll}
    \dsp \frac{Q_{n}(1)}{R_{n}(1)} & \dsp = (1-M_0)\,  \frac{\dsp\sum_{k=1}^{n} \frac{M_{1,0}M_{1,1}\cdots M_{1,k-1}}{(1-M_{1,0})(1-M_{1,1})\cdots (1-M_{1,k-1})}}{\dsp 1 + \sum_{k=1}^{n} \frac{M_{1,0}M_{1,1}\cdots M_{1,k-1}}{(1-M_{1,0})(1-M_{1,1})\cdots (1-M_{1,k-1})}}    , \quad n \geq 1.
  \end{array}
\]
Therefore, Wall's  characterization  for the maximal parameter sequence (see \cite[Thm. 6.2 ]{Chihara-book}) gives us the following lemma.

\begin{lema} \label{Lemma-Convergence-at-one}
\[
      d_1 = \frac{Q_{1}(1)}{R_{1}(1)}  < \frac{Q_{2}(1)}{R_{2}(1)} < \cdots < \frac{Q_{n-1}(1)}{R_{n-1}(1)} < \frac{Q_{n}(1)}{R_{n}(1)} < (1-M_0)\,
\]
and
\[
  \begin{array} {ll}
    \dsp \lim_{n \to \infty}\frac{Q_{n}(1)}{R_{n}(1)}  = (1-M_0).
  \end{array}
\]

\end{lema}

Now we look at the series expansions of the rational functions $Q_{n}(z)/R_{n}(z)$. From the recurrence formula for $\{R_{n}\}$, if $R_{n}(z) = \sum_{j=0}^{n} r_{n,j}\,z^j$, then one can verify that
\begin{equation} \label{Reciprocal-of-R_m}
   R_{n}(z) = \sum_{j=0}^{n} r_{n,j}\,z^j = \sum_{j=0}^{n} \overline{r}_{n,n-j}\,z^j = R_{n}^{\ast}(z), \quad n \geq 0
\end{equation}
and, in particular,
\[ r_{0,0} = 1 \quad \mbox{and} \quad r_{n,n} = \overline{r}_{n,0} = \prod_{k=1}^{n} (1 + i c_k), \ \ n \geq 1.
\]
Here, $R_{n}^{\ast}(z) = z^n \overline{R_{n}(1/\overline{z})}$ is the reciprocal polynomial of $R_{n}(z)$.  With the property (\ref{Reciprocal-of-R_m}) the polynomial $R_n$ can be called a self-inversive polynomial. Likewise, the $n-1$ degree  polynomial $Q_{n}$ also satisfies the self inversive property
\[
      Q_{n}^{\ast}(z) = z^{n-1} \overline{Q_{n}(1/\overline{z})} = Q_{n}(z), \quad n \geq 1.
\]

Applying the respective three term recurrence formulas in %
\[
   U_{n}(z) = \left| \begin{array}{cc}
                     Q_{n}(z) & R_{n}(z) \\
                     Q_{n-1}(z) & R_{n-1}(z)
                     \end{array} \right|
   = Q_{n}(z) R_{n-1}(z) - Q_{n-1}(z) R_{n}(z),  \quad n \geq 1,
\]
there follows \ $U_1(z) = 2d_1$  and
\begin{equation} \label{Determinant-formula}
    U_{n+1}(z) = 4 d_{n+1} z U_{n}(z) = 2^{2n+1} d_1 d_2 \cdots d_{n+1} z^{n}, \quad n \geq 1.
\end{equation}
Such formulas in the literatures of continued fractions and orthogonal polynomials are known as determinant formulas.

Hence, considering the series expansions in terms of the origin and infinity,
\begin{equation} \label{Correspondence-1}
    \frac{Q_{n}(z)}{R_{n}(z)} - \frac{Q_{n-1}(z)}{R_{n-1}(z)} = \left\{
      \begin{array}{l}
        \dsp \frac{\overline{\gamma}_{n-1}}{\overline{r}_{n-1,n-1}}\, z^{n-1} + O\big(z^{n}\big), \\[3ex]
        \dsp \frac{\gamma_{n-1}}{r_{n-1,n-1}}\, \frac{1}{z^{n}} + O\big((1/z)^{n+1}\big),
      \end{array}
      \right.
      n \geq 1,
\end{equation}
where \  $ \dsp \gamma_{n-1} = \frac{2^{2n-1} d_1 d_2 \cdots d_{n}}{r_{n,n}}$, $n \geq 1$. 
This means there exist formal series expansions $E_0(z)$ and $E_{\infty}(z)$, respectively   about the origin and about infinity, such that
\begin{equation} \label{origin-correspondence-1}
    E_{0}(z) - \frac{Q_{n}(z)}{R_{n}(z)} = \frac{\overline{\gamma}_{n}}{\overline{r}_{n,n}}\, z^{n} + O\big(z^{n+1}\big), \quad n \geq 0 \qquad  \quad
\end{equation}
and
\begin{equation} \label{infinity-correspondence-1}
    E_{\infty}(z) - \frac{Q_{n}(z)}{R_{n}(z)} = \frac{\gamma_{n}}{r_{n,n}}\, \frac{1}{z^{n+1}} + O\big((1/z)^{n+2}\big), \quad n \geq 0.
\end{equation}
For reasons that will become clearer later we set
\[
  E_{0}(z) = -\nu_{1} - \nu_{2} z - \nu_{3}z^2 - \nu_{4}z^3 - \ldots \
\]
and
\[
  E_{\infty}(z) = \frac{\nu_{0}}{z} + \frac{\nu_{-1}}{z^2} + \frac{\nu_{-2}}{z^3} + \frac{\nu_{-3}}{z^4} +  \ldots \ .
\]
From (\ref{infinity-correspondence-1}), since  $R_{n}$ and $Q_{n}$ are self inversive, we have
\[
    z^{-1}\overline{E_{\infty}(1/\overline{z})} - \frac{Q_{n}(z)}{R_{n}(z)} = \frac{\overline{\gamma}_{n}}{\overline{r}_{n,n}}\, z^{n} + O\big(z^{n+1}\big),  \quad n \geq 0.
\]
Comparing this with (\ref{origin-correspondence-1}) we then have the symmetry property
\[
    z^{-1}\overline{E_{\infty}(1/\overline{z})}  = E_{0}(z)
\]
and
\[
      \nu_{j} = - \overline{\nu}_{-j+1}, \quad j=1,2, \ldots \ .
\]
Considering the following systems of equations
\[
  \begin{array}{lllllllcl}
  -\nu_{1}r_{n,0}   & & & &  -\ q_{n,0} & & & = & 0, \\[2ex]
  \ \ \ \ \vdots    & \ddots &  &  & & \ddots & & \vdots\\[2ex]
  -\nu_{n}r_{n,0} & \cdots & -\ \nu_{1}r_{n,n-1} &  &  & & -\ q_{n,n-1}& = & 0, \\[2ex]
  -\nu_{n+1}r_{n,0}  & \cdots & -\ \nu_{2}r_{n,n-1} & - \ \nu_{1}r_{n,n}  & & & & = & \overline{\gamma}_{n},
  \end{array}
\]
and
\[
  \begin{array}{lllllllcl}
  \nu_{0}r_{n,0} &+ \ \nu_{-1}r_{n,1} & \cdots  & + \ \nu_{-n}r_{n,n} & & & &  = & \gamma_{n}, \\[2ex]
   & \ \ \ \nu_{0}r_{n,1} & \cdots  & +\ \nu_{-n+1}r_{n,n} & - \ q_{n,0} &  & &  = & 0, \\[2ex]
   &  & \ddots  & \ \ \ \ \ \vdots  &  & \ddots &  &  \vdots & \\[2ex]
   &  & & \ \ \ \nu_{0}r_{n,n} &  &  &-\ q_{n,n-1} &  = & 0
  \end{array}
\]
in  the coefficients of  $Q_n(z)= \sum_{j=0}^{n-1} q_{n,j}z^{j}$ and $R_{n}(z) = \sum_{j=0}^{n} r_{n,j}\,z^j$, which follow from (\ref{origin-correspondence-1}) and (\ref{infinity-correspondence-1}), we have
\[
   \gamma_n = (-1)^{n} \frac{H_{n+1}^{(-n)}}{\overline{H}_{n}^{(-n+1)}}\, r_{n,n}, \quad n \geq 1,
\]
where
\[
 H_{n}^{(m)} =
  \left|
   \begin{array}{cccc}
   \nu_m & \nu_{m+1} & \cdots  &  \nu_{m+n-1} \\
   \nu_{m+1} & \nu_{m+2} &  \cdots &  \nu_{m+n} \\
   \vdots  & \vdots  &    &  \vdots  \\
   \nu_{m+n-1} & \nu_{m+n}  & \cdots &  \nu_{m+2n-2}
   \end{array}
  \right|,
\]
are the Hankel determinants associated with the double sequence $\{\nu_n\}_{n=-\infty}^{\infty}$.  Moreover, the following lemma can also be stated.

\begin{lema} \label{Lemma-Orto-Rn}
Let $\{R_n\}$ and $\{Q_n\}$ be the sequences of polynomials obtained from  the real sequence  $\{c_n\}$, the positive chain sequence $\{d_{n}\}$ and the three term recurrence formulas $(\ref{Eq-TTRR-Rn})$ and $(\ref{Eq-TTRR-Qn})$. Then the rational functions $Q_n/R_n$ satisfy the correspondence properties given by $(\ref{origin-correspondence-1})$ and $(\ref{infinity-correspondence-1})$. Moreover, if the moment functional $\mathcal{N}$ is such that
\[
    \mathcal{N}[\zeta^{-n}] = \nu_{n}, \quad n=0, \pm1,\pm2, \ldots \ ,
\]
then the polynomials $R_n$  satisfy
\[
    \mathcal{N}[\zeta^{-n+j} R_{n}(\zeta)] = \left\{
      \begin{array}{cl}
        -\overline{\gamma}_{n}, &  j=-1,\\[1ex]
        0, & j=0,1,\ldots, n-1, \\[1ex]
        \gamma_{n}, & j=n
      \end{array} \right. \quad n \geq 1,
\]
where \quad $\dsp \gamma_0 = \nu_0 = \frac{2d_1}{1+ic_1}$\quad and \quad $\dsp \gamma_n = \frac{4d_{n+1} }{(1+ic_{n+1})} \gamma_{n-1}$,\ $n\geq 1$. \\

\end{lema}

In terms of $\gamma_{n} = \mathcal{N}[R_{n}(\zeta)]$, $n \geq 0$, the elements of  $\{c_n\}$ and $\{d_{n}\}$ satisfy
\[
    \frac{1}{d_1} = \frac{1}{\gamma_{0}} + \frac{1}{\overline{\gamma}_{0}}, \qquad i \frac{c_1}{d_1} = \frac{1}{\gamma_{0}} - \frac{1}{\overline{\gamma}_{0}}
\]
and
\[
       \frac{1}{2d_{n+1}} =  \frac{\gamma_{n-1}}{\gamma_{n}} + \frac{\overline{\gamma}_{n-1}}{\overline{\gamma}_{n}},  \qquad  i \frac{c_{n+1}}{2d_{n+1}} =  \frac{\gamma_{n-1}}{\gamma_{n}} - \frac{\overline{\gamma}_{n-1}}{\overline{\gamma}_{n}},  \qquad n \geq 1.
\]

Observe also that one can write
\[ \mathcal{N}\big[\frac{1}{z-\zeta}\big] = \left\{
      \begin{array}{ll}
        E_{0}(z), &  \mbox{for} \ z \sim 0,\\[1ex]
        E_{\infty}(z), &  \mbox{for} \ z \sim \infty.
      \end{array} \right.
\]
%

%%%%%%%%%%%%%%%%%%%%%%%%%%%%%%%%%%%%%%%%%%%%%%%%%%%%%%%%%%%%%%%
\setcounter{equation}{0}
\section{Associated  moments and measure on the unit circle} \label{Sec-Recovering-Szego}

Given the sequence of polynomials $\{R_{n}\}$ and $\{Q_{n}\}$, as defined in sections \ref{Sec-Intro} and \ref{Sec-CorrAsymProps}, let the sequence of polynomials  $\{A_{n}\}$ and $\{B_{n}\}$ be given by
\begin{equation} \label{Szego-RationalFunction}
     A_{n}(z) = R_{n}(z) - Q_{n}(z) \quad \mbox{and} \quad B_{n}(z) = (z-1)R_{n}(z),  \quad n \geq 0.
\end{equation}
It is easily verified that
\[
   \frac{A_{n+1}(z)}{B_{n+1}(z)} - \frac{A_{n}(z)}{B_{n}(z)}= -\frac{1}{z-1} \left[\frac{Q_{n+1}(z)}{R_{n+1}(z)} - \frac{Q_{n}(z)}{R_{n}(z)}  \right], \quad n \geq 0.
\]
Hence, from (\ref{Correspondence-1})
\begin{equation} \label{Correspondence-2}
    \frac{A_{n+1}(z)}{B_{n+1}(z)} - \frac{A_{n}(z)}{B_{n}(z)} = \left\{
      \begin{array}{l}
        \dsp \ \frac{\overline{\gamma}_{n}}{\overline{r}_{n,n}}\, z^{n} + O\big(z^{n+1}\big), \\[3ex]
        \dsp -\frac{\gamma_{n}}{r_{n,n}}\, \frac{1}{z^{n+2}} + O\big((1/z)^{n+3}\big),
      \end{array}
      \right.
      n \geq 0.
\end{equation}
Thus, we can state the following lemma.

\begin{theo}  \label{Thm-MainCorrespondence}
Associated with the real sequence $\{c_n\}$ and the positive chain sequence $\{d_n\}$ there exists a nontrivial probability measure $\mu$ on the unit circle. If  $M_0 > 0$, where $\{M_n\}$ is the maximal parameter sequence of $\{d_n\}$, then $\mu$ has a pure point of mass $M_0$ at $z=1$. Let $\mathcal{N}$ be the moment functional associated with $\{c_n\}$ and $\{d_n\}$ as given by Lemma \ref{Lemma-Orto-Rn}. Then 
\[
     \mathcal{N}[\zeta^{-n}] = \int_{\mathcal{C}} \zeta^{-n} (1-\zeta) d \mu(\zeta), \quad n = 0, \pm1, \pm2, \ldots \ . 
\]

\end{theo}

\noindent {\bf Proof}. From (\ref{Correspondence-2}) there exist series expansions $F_0(z)$ and $F_{\infty}(z)$ such that
\begin{equation} \label{origin-correspondence-2}
    F_{0}(z) - \frac{A_{n}(z)}{B_{n}(z)} = \frac{\overline{\gamma}_{n}}{\overline{r}_{n,n}}\, z^{n} + O\big(z^{n+1}\big), \quad n \geq 0 \qquad  \quad
\end{equation}
and
\begin{equation} \label{infinity-correspondence-2}
    F_{\infty}(z) - \frac{A_{n}(z)}{B_{n}(z)} = - \frac{\gamma_{n}}{r_{n,n}}\, \frac{1}{z^{n+2}} + O\big((1/z)^{n+3}\big), \quad n \geq 0.
\end{equation}
Setting
\[
  F_{0}(z) = -\mu_{1} - \mu_{2} z - \mu_{3} z^2 - \mu_{4}z^3 - \ldots \ ,
\]
and
\[
  F_{\infty}(z) = \frac{\mu_{0}}{z} + \frac{\mu_{-1}}{z^2} + \frac{\mu_{-2}}{z^3} + \frac{\mu_{-3}}{z^4} +  \ldots \ ,
\]
we then have from (\ref{origin-correspondence-1}), (\ref{infinity-correspondence-1}), (\ref{Szego-RationalFunction}) (\ref{origin-correspondence-2}) and (\ref{infinity-correspondence-2}) that the numbers $\mu_k$ satisfy
\begin{equation} \label{Eq-Moment-Relations}
  \begin{array}{l}
      \mu_{n} = 1 + \sum_{j=1}^{n} \nu_{j},   \\[2ex]
     \mu_{-n} = 1 - \sum_{j=1}^{n} \nu_{-j+1},
  \end{array} n = 1, 2, 3, \ldots \ ,
\end{equation}
with $\mu_0 = 1$. Since $\nu_{j} = - \overline{\nu}_{-j+1}$, $j \geq 1$, there hold
\[
  \mu_{-k} = \overline{\mu}_{k}, \quad k \geq 1.
\]
If we define the moment functional $\mathcal{M}$ by
\[
     \mathcal{M}[\zeta^{-k}] = \mu_{k}, \quad k=0, \pm1, \pm2, \ldots \ ,
\]
then
\begin{equation} \label{Eq-Functional-Relations-1}
    \mathcal{M}[\zeta^{-k}] =  1 - \mathcal{N}\big[\frac{1-\zeta^{-k}}{1-\zeta}\big], \quad k=0, \pm1, \pm2, \ldots \ .
\end{equation}
Since $\nu_{-k} = \mu_{-k} - \mu_{-k-1}$, $k=0, \pm1, \pm2, \ldots $, which follows from (\ref{Eq-Moment-Relations}), the moment functionals $\mathcal{M}$ and $\mathcal{N}$ also satisfy
\begin{equation} \label{Eq-Functional-Relations-2}
     \mathcal{N}[\zeta^{k}] =  \mathcal{M}[\zeta^{k}(1-\zeta)] , \quad  k=0, \pm1, \pm2, \ldots \ .
\end{equation}

Now considering the partial decomposition of $A_{n}(z)/B_{n}(z)$, we have
\begin{equation} \label{Eq-Partial-Decomp}
   \frac{A_n(z)}{B_n(z)} = \frac{R_{n}(z) - Q_{n}(z)}{(z-1)R_{n}(z)} = \frac{\lambda_{n,0}}{z-1} + \sum_{j=1}^{n} \frac{\lambda_{n,j}}{z-z_{n,j}},
\end{equation}
where $z_{n,j}$, $j=1,2, \ldots, n$, are the zeros of $R_n(z)$,
\[
    \lambda_{n,0} = 1 - \frac{Q_{n}(1)}{R_{n}(1)}
\]
and
\[
   \lambda_{n,j} = \frac{Q_{n}(z_{n,j})}{(1-z_{n,j})R_{n}^{\prime}(z_{n,j})}, \quad j =1, 2, \ldots n.
\]
Clearly, from Lemma \ref{Lemma-Convergence-at-one},
\[
     1 - d_1 = \lambda_{1,0} > \lambda_{2,0} > \cdots > \lambda_{n,0} > \lambda_{n+1,0} > \cdots
\]
and
\[
   \lim_{n \to \infty} \lambda_{n,0} = M_0.
\]
Furthermore, since we can write
\[
   \lambda_{n,j} = \frac{U_{n}(z_{n,j})}{(1-z_{n,j})V_{n}(z_{n,j})}, \quad j =1, 2, \ldots n,
\]
where $V_n(z)$ and $U_n(z)$ are respectively given by (\ref{Wronskian-for-General-Rn}) and (\ref{Determinant-formula}),
we have
\[
     \lambda_{n,j} = \frac{4d_1 d_2 \cdots d_{n}}{W_n(x_{n,j})} \frac{z_{n,j}}{(z_{n,j} - 1)(1 - z_{n,j})} > 0, \quad j =1,2, \ldots, n.
\]
Here, $W_n(x)$ are the Wronskians defined in (\ref{Wronskian-for-General-Gn}) and, with $z_{n,j} = e^{i\theta_{n,j}}$,
\[
     \frac{z_{n,j}}{(z_{n,j} - 1)(1 - z_{n,j})} = \frac{1}{4 \sin^2(\theta_{n,j}/2)}.
\]
In addition to the positiveness of the elements $\lambda_{n,j}$, $j=0,1, 2, \ldots, n$, by considering the limit of  $zA_{n}(z)/B_{n}(z)$, as $z \to \infty$, we also have
\[
     \sum_{j=0}^{n} \lambda_{n,j} = 1.
\]

Now if the step functions $\psi_{n}(e^{i\theta})$, $n \geq 1$, are defined on $[0, 2\pi]$ by
\[
        \psi_n(e^{i\theta}) = \left\{
        \begin{array}{ll}
          0, & \theta = 0, \\
          \lambda_{n,0}, & 0 < \theta \leq \theta_{n,1}, \\
          \sum_{j=0}^{k} \lambda_{n,j}^{(t)}, & \theta_{n,k} < \theta \leq \theta_{n,k+1},\quad k=1,2, \ldots, n-1, \\
          1, & \theta_{n,n} < \theta \leq 2\pi.
        \end{array} \right.
\]
then from the definition of the Riemann-Stieltjes integrals
\[
     \frac{A_n(z)}{B_n(z)} = \int_{\mathcal{C}} \frac{1}{z - \zeta}\, d \psi_{n}(\zeta), \quad n \geq 1.
\]
Hence, by the application of the Helley selection theorem (see \cite{JoNjTh-1989}) there exists a subsequence $\{n_j\}$ such that $\psi_{n_j}(e^{i\theta})$ converges to a bonded non-decreasing function, say $\mu(e^{i\theta})$, in $[0, 2\pi]$.

From (\ref{origin-correspondence-2}) and (\ref{infinity-correspondence-2}), since
\[
   \int_{\mathcal{C}} d \psi_{n}(\zeta) = 1 \quad \mbox{and} \quad \int_{\mathcal{C}} \zeta^{k} d \psi_{n}(\zeta) = \mu_{-k}, \ \ k =\pm1,\pm2, \ldots, \pm n,
\]
we also have that
\[
   \int_{\mathcal{C}} d \mu(\zeta) = 1 = \mathcal{M}[1] \quad \mbox{and} \quad \int_{\mathcal{C}} \zeta^{k} d \mu(\zeta) = \mu_{-k}= \mathcal{M}[\zeta^k], \ \ k =\pm 1,\pm 2, \ldots\, .
\]
Now, $\mu$ is the only probability measure  that satisfies the above relations follows from known results on the moment problem on the unit circle. The measure has a jump $M_0$ at $z=1$ is also confirmed by  $\lim_{n \to \infty} \lambda_{n,0} = M_0$.  \hfill \qed

%%%%%%%%%%%%%%%%%%%%%%%%%%%%%%%%%%%%%%%%%%%%%%%%%%%%%%%%%%%%%%%
\setcounter{equation}{0}
\section{Further properties of $R_n(z)$ and the associated OPUC} \label{Sec-ConcRem}

With the probability measure $\mu$ obtained in the previous section we then have from Lemma \ref{Lemma-Orto-Rn} and Theorem \ref{Thm-MainCorrespondence} that 
\[
   \nu_k = \mathcal{N}[\zeta^{-k}] = \int_{\mathcal{C}}\zeta^{-k} (1-\zeta) d \mu(\zeta), \quad k = 0, \pm1, \pm2, \ldots
\]
and for $n \geq 1$,
\begin{equation} \label{Eq-Orthog-Rn}
    \int_{\mathcal{C}}\zeta^{-n+k} R_n(\zeta)(1-\zeta) d \mu(\zeta) = 0, \quad 0 \leq k \leq n-1.
\end{equation}
The following lemma provides information about the values of the integrals \linebreak  $\int_{\mathcal{C}}R_n(\zeta) d \mu(\zeta)$.

\begin{lema} \label{Lemma-Another-Param-PCS}  Let \ 
\(
     \dsp \widehat{\gamma}_n = \int_{\mathcal{C}}R_n(\zeta) d \mu(\zeta), \ n \geq 0. \ 
\)
Then
\begin{equation} \label{Eq-Another-Param-PCS}
     \widehat{\gamma}_0 = 1 \quad \mbox{and} \ \quad \widehat{\gamma}_{n} = 2(1-m_n)\,\widehat{\gamma}_{n-1}, \ n \geq 1,
\end{equation}
where $\{m_n\}_{n=0}^{\infty}$ is the minimal parameter sequence of the positive chain sequence $\{d_n\}$.
\end{lema}

\noindent {\bf Proof}.  First we observe from  (\ref{Eq-Orthog-Rn}) that
\begin{equation} \label{Eq-Param-RR}
  \widehat{\gamma}_n = \int_{\mathcal{C}}\zeta^{-k}R_n(\zeta) d \mu(\zeta), \quad k=0,1, \ldots, n, \quad n \geq 1.
\end{equation}
Now by direct evaluations $\widehat{\gamma}_0 = \int_{\mathcal{C}}d \mu(\zeta) = 1$ and
\[
    \widehat{\gamma}_1 = (1+ic_1)\mu_{-1} + (1-ic_1).
\]
Thus, from $\mu_{-1} = 1 - \nu_0$ and $\nu_0 = 2d_1/(1+ic_1)$, there follows
\[
     \widehat{\gamma}_1 = 2(1-d_1) = 2(1-m_1),
\]
proving (\ref{Eq-Another-Param-PCS}) for $n =1$.

Now from the three term recurrence formula (\ref{Eq-TTRR-Rn}), we have
\[
  \begin{array}{ll}
   \int_{\mathcal{C}}\zeta^{-1}R_{n+1}(\zeta) d \mu(\zeta)  = & \int_{\mathcal{C}}\zeta^{-1}(\zeta+1)R_{n}(\zeta) d \mu(\zeta) \\[2ex]
   & \qquad + \ ic_{n+1} \int_{\mathcal{C}}\zeta^{-1}(\zeta-1)R_{n}(\zeta) d \mu(\zeta) \\[2ex]
   & \qquad \qquad \qquad \qquad - \ 4\,d_{n+1}\int_{\mathcal{C}}R_{n-1}(\zeta) d \mu(\zeta),
  \end{array}
\]
for $n \geq 1$. Using (\ref{Eq-Orthog-Rn}) and (\ref{Eq-Param-RR}) we then have
\[
    \widehat{\gamma}_{n+1} = 2 \widehat{\gamma}_{n} - 4 d_{n+1} \widehat{\gamma}_{n-1}, \quad n \geq 1.
\]
From this we have
\[
    d_{n+1} = d_{n,1} = \frac{\widehat{\gamma}_{n}}{2\widehat{\gamma}_{n-1}}\big(1 - \frac{\widehat{\gamma}_{n+1}}{2\widehat{\gamma}_{n}}\big), \quad n \geq 1.
\]
Since the minimal parameter sequence of the positive chain sequence $\{d_{n}\}$ is also a parameter sequence of the positive chain sequence  $\{d_{n,1}\}$, we obtain  from
\[ \frac{\widehat{\gamma}_{1}}{2\widehat{\gamma}_{0}} = (1 - m_1) , \]
that
\[
    \frac{\widehat{\gamma}_{n}}{2\widehat{\gamma}_{n-1}} = (1-m_n), \quad n \geq 1,
\]
which completes the proof of the lemma. \hfill \qed

As an immediate consequence of the above results we have
\[
    \int_{\mathcal{C}}{\overline{\zeta^{k}}}[R_{n}(\zeta) - 2(1-m_n)R_{n-1}(\zeta)] d \mu(\zeta) = \int_{\mathcal{C}}\zeta^{-k}[R_{n}(\zeta) - 2(1-m_n)R_{n-1}(\zeta)] d \mu(\zeta) = 0,
\]
for $k=0,1,\ldots, n-1$.  Since $R_{n}(z) - 2(1-m_n)R_{n-1}(z)$ is a polynomial of exact degree $n$ with leading coefficient $r_{n,n} = \prod_{k=1}^{n}(1+ic_k)$, we can state the following.

\begin{theo} \label{Thm-Main} 
If the sequence of polynomials $\{S_n\}$ is such that
\[
   S_0(z) = 1 \quad \mbox{and} \quad S_n(z) \prod_{k=1}^{n}(1+ic_k) = R_{n}(z) - \,2(1-m_n)R_{n-1}(z), \quad n \geq 1,
\]
then $\{S_n\}$ is the  sequence of monic OPUC with respect to the measure $\mu$.

\end{theo}

From the above theorem, together with the formula for $R_n(0)$ given after  (\ref{Reciprocal-of-R_m}), the associated Verblunsky coefficients $\alpha_{n-1} = - \overline{S_n(0)}$, $n \geq 1$, are
\[
      \alpha_{n-1} = \frac{1-2m_n-ic_n}{1+ic_n} \prod_{k=1}^{n}\frac{1 + ic_k}{1 - i c_k}, \quad n \geq 1.
\]

%%%%%%%%%%%%%%%%%%%%%%%%%%%%%%%%%%%%%%%%%%%%%%%%%%%%%%%%%%%%%%%
\setcounter{equation}{0}
\section{An example} \label{Sec-ConcRem}

We now analyze the results obtained so far with the use of the following example. \\

Let the real sequences $\{c_n\}$ and $\{d_n\}$ be  given by
\begin{equation} \label{Eq-Special-Example}
 \begin{array}l
  \dsp c_n = \frac{\eta}{\lambda+n}, \ \ n \geq 1, \\[2ex]
  \dsp d_1 = d_1(t) = \frac{1}{2} \frac{2\lambda+1}{\lambda+1} (1-t), \quad d_{n+1} = \frac{1}{4} \frac{n(2\lambda+n+1)}{(\lambda+n)(\lambda+n+1)}, \ \ n\geq 1,
 \end{array}
\end{equation}
where \ $\lambda, \eta \in \mathbb{R}$, \ $\lambda > -1/2$ \ and \ $0 \leq t < 1$.

Observe that, as verified in \cite{Costa-Felix-Ranga-2013}, $\{d_n\}_{n=1}^{\infty}$ is a positive chain sequence with its maximal parameter sequence $\{M_n^{(t)}\}_{n=0}^{\infty}$ given by
\begin{equation} \label{Eq-Special-ParamSeq}
   M_0^{(t)} = t, \quad M_n^{(t)}  = \frac{1}{2} \frac{2\lambda+n}{\lambda+n}, \ \ n \geq 1.
\end{equation}

It was first shown in \cite{Ranga-2010} (see also \cite{Costa-Felix-Ranga-2013}) that the polynomials $R_n$ obtained from the above sequences $\{c_n\}$ and $\{d_n\}$, together with the recurrence formula (\ref{Eq-TTRR-Rn}), are
\[
      R_n(z) = \frac{(2\lambda+2)_n}{(\lambda+1)_n} \, _2F_1(-n,b+1;\,b+\overline{b}+2;\,1-z),  \quad n \geq 1,
\]
where  $b = \lambda + i \eta$. This actually follows from the contiguous relation
\[
  \begin{array}{ll}
   (c-a)\,_2F_1(a-1,b;\,c;\,z) & \dsp = \big(c-2a - (b-a)z\big)\, _2F_1(a,b;\,c;\,z) \\[1ex]
               &  \dsp \qquad  \qquad  \qquad \qquad   +\ a (1-z)\,  _2F_1(a+1,b;\,c;\,z),
  \end{array}
\]
of Gauss (see  \cite[Eq. (2.5.16)]{Andrews-Book}), by letting $a = -n$, $b = b+1$, $c = b+\overline{b}+2$ and $z = 1-z$.

Now consider the rational functions $Q_n/R_n$, $n \geq 1$, given by the continued fraction expression
\begin{equation} \label{Special-CF-Expansion}
  \begin{array} {l}
    \dsp \frac{Q_{n}(z)}{R_{n}(z)}  = \cfr{2d_1}{(1+ic_1)z+(1-ic_1)}-\cfr{4d_2z}{(1+ic_2)z+(1-ic_2)} - \cdots   \\[2ex]
     \dsp \qquad \qquad \qquad \qquad \qquad \qquad \qquad \qquad \qquad \qquad \qquad \qquad  \cdots - \cfr{4d_{n}z}{(1+ic_n)z+(1-ic_n)} \,,
  \end{array}
\end{equation}
Clearly, by Lemma \ref{Lemma-Convergence-at-one}, the increasing sequence $\{Q_n(1)/R_n(1)\}$ satisfy
\[
  \begin{array} {ll}
    \dsp \lim_{n \to \infty}\frac{Q_{n}(1)}{R_{n}(1)}  = 1-t.
  \end{array}
\]
In order to obtain the power series expansions associated with the rational functions $Q_n(z)/R_n(z)$, we write the continued fraction expression (\ref{Special-CF-Expansion}) in the equivalent form
\begin{equation} \label{New-CF-Expansion}
  \begin{array} {l}
    \dsp \frac{(1-ic_1)Q_{n}(z)}{2d_1\,R_{n}(z)}  = \cfr{1}{1+b_1z}-\cfr{a_2z}{1 +b_2z}   -  \cdots - \cfr{a_nz}{1 + b_n z} \, ,
  \end{array}
\end{equation}
where
\[
    b_n = \frac{1+ic_n}{1-ic_n} = \frac{b+n}{\overline{b}+n}, \quad a_{n+1} = \frac{4d_{n+1}}{(1-ic_n)(1-ic_{n+1})} = \frac{n(b+\overline{b}+n+1)}{(\overline{b}+n)(\overline{b}+n+1)}, \quad n \geq 1.
\]
We thus can obtain  the power series expansion about the origin of the rational function on the left hand side of (\ref{New-CF-Expansion}) from
\[
 \Omega_{0}(b;\, z) = \cfr{1}{1+b_1z}-\cfr{a_2z}{1 +b_2z}   -  \cdots - \cfr{a_{n-1}z}{1 + b_{n-1} z} - \cfr{a_{n}z}{1 + b_{n} z - a_{n+1}z\, \Omega_{n}(b;\, z)} \, ,
\]
where $\dsp \Omega_{n}(b;\, z) = \frac{\,_2F_1(n+1,-b;\,\overline{b}+n+2;\,z)}{\,_2F_1(n,-b;\,\overline{b}+n+1;\,z)}$.

The above relation follows from the contiguous relation
\[
  \begin{array}{ll}
   _2F_1(a,b;c;z) & \dsp = \Big( 1 + \frac{a-b+1}{c}z\Big)\,_2F_1(a+1,b;c+1;z) \\[2ex]
               &  \dsp \qquad  \qquad  \qquad \qquad \qquad  - \frac{(a+1)(c-b+1)}{c(c+1)} z\,_2F_1(a+2,b;c+2;z)
  \end{array}
\]
of Gauss (see \cite[Eq. (2.5.3)]{Andrews-Book}), by substituting $a$, $b$ and $c$ with $n$, $-b$ and  $\overline{b}+n+1$, respectively.
Since $\Omega_{0}(b;\, z) = \,_2F_1(1,-b;\,\overline{b}+2;\,z)$, we have
\[
    \frac{Q_{n}(z)}{R_{n}(z)} \sim \frac{2d_1}{1-ic_1}\,_2F_1(1,-b;\,\overline{b}+2;\,z)  = -\nu_{1} - \nu_{2} z - \nu_{3}z^2 - \ldots \, ,
\]
from which
\begin{equation} \label{Special-Modified-Moments}
     \nu_n = \frac{b+\overline{b}+1}{b+1} \frac{(-b-1)_{n}}{(\overline{b}+1)_n}(1-t), \quad n \geq 1.
\end{equation}
Since $\nu_{n} = - \overline{\nu}_{-n+1}$, with the convention $(a)_n = \Gamma(a+n)/\Gamma(a)$ for all integer values of $n$, the   above expression for $\nu_n$ is also valid for $n \leq 0$.

\begin{figure}[b]
%\tabcapfont
\centerline{
\begin{tabular}{c}
\includegraphics[width=5cm,height=5cm]{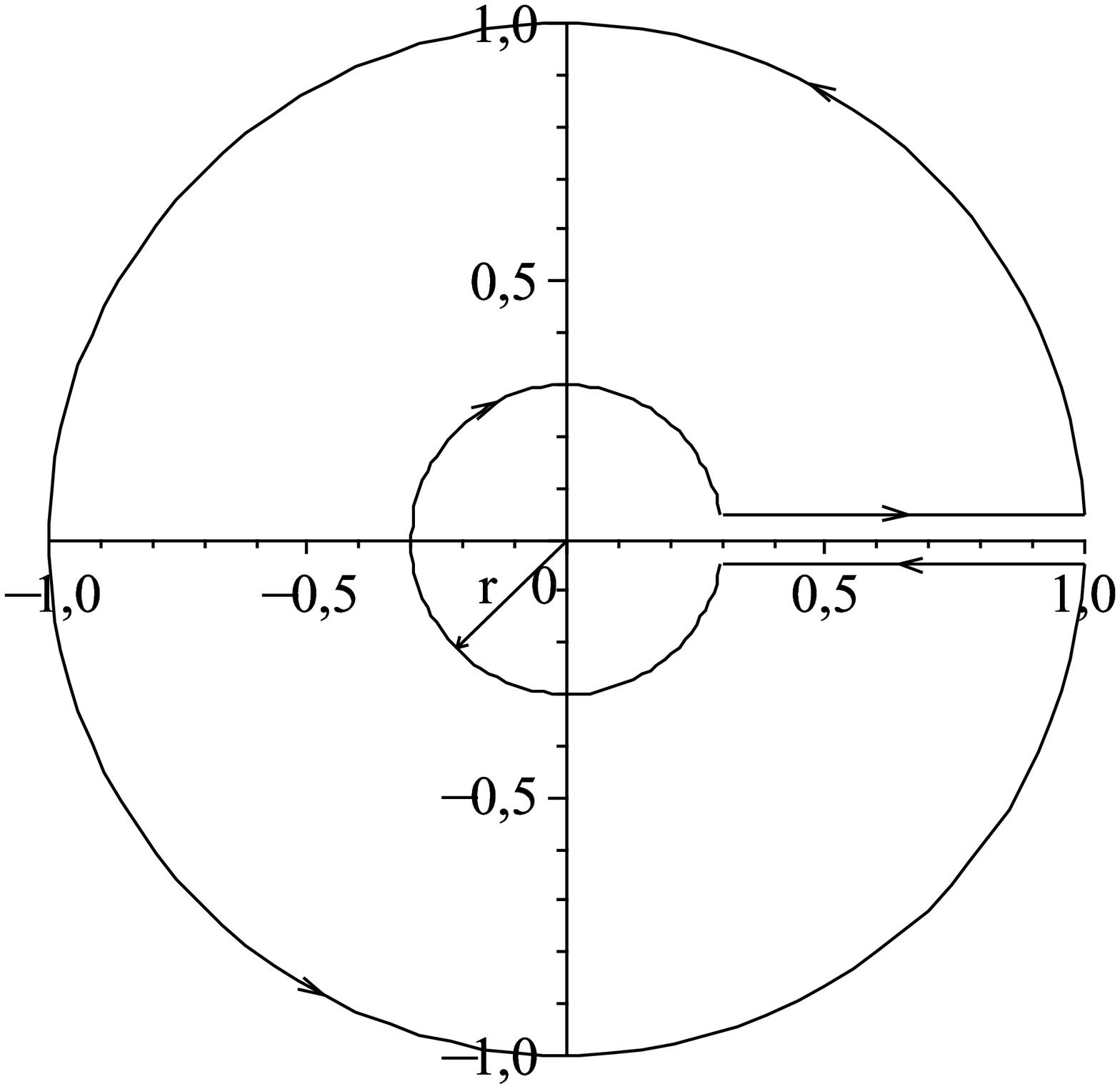}
 \end{tabular}}
\caption{Contour $\Lambda$ } \label{fig1}
\end{figure}

Now we show that
\begin{equation} \label{Integral-Representation-for-SMM}
  \begin{array}{l}
    \dsp \nu_n = \mathcal{N}[\zeta^{-n}]  =
      \int_{\mathcal{C}}\zeta^{-n} \,d \psi(t,b;\,\zeta), \quad n = 0, \pm1,\pm2, \ldots \ ,  \\[2ex]
  \end{array}
\end{equation}
where \quad
\(
  \dsp d \psi(t,b;\,z) = \frac{i\,|\Gamma(b+1)|^2(1-t)}{2\pi\, \Gamma(b+\overline{b}+1)} (-z)^{-\overline{b}-1} (1-z)^{b+\overline{b}+1} dz .
\)

\noindent Here, the branch cuts in $(-z)^{-\overline{b}} = (e^{-i\pi}z)^{-\overline{b}}$ and $(1-z)^{b+\overline{b}} = \big(e^{-i\pi}(z-1)\big)^{b+\overline{b}}$ are along the positive real axis.

Before verifying (\ref{Integral-Representation-for-SMM}), we observe that (\ref{Integral-Representation-for-SMM}) can also be written in the equivalent form
\[
    \nu_n = \frac{-i\,2^{b+\overline{b}+1}|\Gamma(b+1)|^2(1-t)}{2\pi\, \Gamma(b+\overline{b}+1)}\int_{0}^{2\pi} e^{-in\theta} e^{-i\theta/2} e^{(\pi-\theta)\mathfrak{Im}(b)} [\sin^{2}(\theta/2)]^{\mathfrak{Re}(b)+1/2} d\theta . \qquad
\]

First we show (\ref{Integral-Representation-for-SMM}) for $\nu_{-n}$ for those values of $n$ such  that $\mathfrak{Re}(n-\overline{b}) >  0$. Since $\int_{\Lambda}(-z)^{n-\overline{b}-1} (1-z)^{b+\overline{b}+1} dz = 0$, where $\Lambda$ is the contour given as in Fig.\,\ref{fig1}, we obtain for $\mathfrak{Re}(n-\overline{b}) >  0$ that

\[
      \nu_{-n} = \frac{-2\sin(\overline{b}\pi)\,|\Gamma(b+1)|^2(1-t)}{2\pi\, \Gamma(b+\overline{b}+1)}
      \int_{0}^{1}x^{n-\overline{b}-1} (1-x)^{b+\overline{b}+1} dx.
\]
Hence,  from   the definitions of the gamma function, the beta function and the Euler's reflection formula, we obtain  the required result.

To obtain the result for other values of $n$, we note that for $\nu_n$ given by (\ref{Integral-Representation-for-SMM}) there hold
\[
     \nu_{n} = \frac{\overline{b}+ 1 + n}{-b-1+n} \nu_{n+1} \quad \mbox{and}\quad   \nu_{n} = - \overline{\nu}_{-n+1}, \quad n = 0, 1,2, \ldots ,
\]
the first of these results follows from integration by parts and other by simple conjugation.

The idea used here to calculate the integral (\ref{Integral-Representation-for-SMM}) is the same employed in, for example,  \cite{Hend-vanRo-1986, Ranga-2010}. In  \cite{Hend-vanRo-1986} the authors consider a  general set of parameters for the exponents of $z$ and $1-z$, but restricting the values of the parameters to be real.

Now from the representation (\ref{Special-Modified-Moments}) for the coefficients $\nu_n$, we have from (\ref{Eq-Moment-Relations}) that
\[
    \mathcal{M}[1] = 1, \quad \overline{\mathcal{M}[\zeta^{n}]} = \mathcal{M}[\zeta^{-n}]  = \mu_n = 1 + \frac{b+\overline{b}+1}{b+1} (1-t)  \sum_{j=1}^{n} \frac{(-b-1)_{j}}{(\overline{b}+1)_j}, \quad n \geq 1.
\]
Hence,
\[
    \overline{\mu}_{-n} =  \mu_n = t +  (1-t)  \frac{(-b)_{n}}{(\overline{b}+1)_n} , \quad n \geq 0 ,
\]
which follows from the interesting summation formula
\[
    1 + \frac{b+\overline{b}+1}{b+1}   \sum_{j=1}^{n} \frac{(-b-1)_{j}}{(\overline{b}+1)_j} = \frac{(-b)_{n}}{(\overline{b}+1)_n} , \quad n \geq 1,
\]
verified easily by induction.

Now from the representation  (\ref{Integral-Representation-for-SMM}) for the coefficients $\nu_n$, we have from  (\ref{Eq-Functional-Relations-1}) that
\[
  \mathcal{M}[\zeta^{-n}]  = 1 - (1-t)\frac{i\,|\Gamma(b+1)|^2}{2\pi\, \Gamma(b+\overline{b}+1)} \int_{\mathcal{C}}(1-\zeta^{-n}) \, (-\zeta)^{-\overline{b}-1} (1-\zeta)^{b+\overline{b}} d\zeta , \quad
\]
for $n = 0, \pm1,\pm2, \ldots \ $. Since we can verify using  integration by parts
\[
  (1-t) \frac{i\,|\Gamma(b+1)|^2}{2\pi\, \Gamma(b+\overline{b}+1)} \int_{\mathcal{C}} \, (-\zeta)^{-\overline{b}-1} (1-\zeta)^{b+\overline{b}} d\zeta = -\frac{\overline{b}+1}{b+\overline{b}+1}\, \nu_{1} = 1 - t ,
\]
we can write

\begin{equation} \label{Special-Positive-Measure}
  \begin{array}{ll}
    \dsp \mathcal{M}[\zeta^{-n}] & \dsp  = \int_{\mathcal{C}}\zeta^{-n} d \mu(t,b;\,\zeta) \\[1ex]
    &\dsp    = t (1^{-n}) + (1-t) \frac{i\,|\Gamma(b+1)|^2}{2\pi\, \Gamma(b+\overline{b}+1)} \int_{\mathcal{C}}\zeta^{-n} \, (-\zeta)^{-\overline{b}-1} (1-\zeta)^{b+\overline{b}} d\zeta ,
  \end{array}
\end{equation}
for $n = 0, \pm1,\pm2, \ldots \ $.  Equivalently, this can also be written as
\begin{equation*}
    \dsp \mu_n = \mathcal{M}[\zeta^{-n}] = t\,e^{in0} + (1-t)\,\frac{2^{b+\overline{b}}|\Gamma(b+1)|^2}{2\pi\, \Gamma(b+\overline{b}+1)}\int_{0}^{2\pi} e^{-in\theta}\, [e^{(\pi-\theta)}]^{\mathfrak{Im}(b)} [\sin^{2}(\theta/2)]^{\mathfrak{Re}(b)} d\theta, \qquad
\end{equation*}
for $n = 0, \pm1,\pm2, \ldots \ $.

By Theorem \ref{Thm-Main}, the monic OPUC associated with the positive measure $\mu(t,b;\,z)$ given by (\ref{Special-Positive-Measure}) are $S_0^{(t)}(z) = 1$ and
\[
  \begin{array}l
      \dsp S_n^{(t)}(z) = \frac{(2\lambda+2)_n}{(b+1)_n} \Big[\, _2F_1(-n,b+1;\,b+\overline{b}+2;\,1-z) \\[1ex]
      \dsp \qquad \qquad \qquad \qquad \qquad  - 2(1-m_n^{(t)}) \frac{(\lambda+n)}{2\lambda+n+1}\, _2F_1(-n+1,b+1;\,b+\overline{b}+2;\,1-z)  \Big], \ \ n \geq 1,
  \end{array}
\]
where $\{m_n^{(t)}\}_{n=0}^{\infty}$, such that
\[
   m_0^{(t)} = 0, \qquad m_n^{(t)} = d_n/(1-m_{n-1}^{(t)}), \ \ n \geq 1,
\]
is the minimal parameter sequence of the positive chain sequence given in  (\ref{Eq-Special-Example}).

When $t = 0$, the minimal $\{m_n^{(t)}\}_{n=0}^{\infty}$ and maximal $\{M_n^{(t)}\}_{n=0}^{\infty}$ parameter sequences of $\{d_n\}$ coincide and, hence,  we obtain from (\ref{Eq-Special-ParamSeq}) that $S_0^{(0)}(z) = 1$ and
\[
  \begin{array}l
      \dsp S_n^{(0)}(z) = \frac{(2\lambda+1)_n}{(b+1)_n} \, _2F_1(-n,b+1;\,b+\overline{b}+1;\,1-z) ,  \ \ n\geq 1,
  \end{array}
\]
these are the orthogonal polynomials on the unit circle (see \cite{Ranga-2010}) with respect to the positive measure $\mu(0,b;\,z)$ given by (\ref{Special-Positive-Measure}).

%%%%%%%%%%%%%%%%%%%%%%%%%%%%%%%%%%%%%%%%%%%%%%%%%%%%%%%%%

\end{document}